\documentstyle[12pt]{article}
\def\one{{1\hskip-4pt 1}}
\def\ptl{\partial}
\def\ad{\mbox{ad}}
\def\DD{\hbox{$I\hskip -4pt D$}}

\def \Z{\hbox{$Z\hskip -5.2pt Z$}}

\def\sZ{\hbox{$\sc Z\hskip -4.2pt Z$}}
\def \Q{\hbox{$Q\hskip -5pt \vrule height 6pt depth 0pt\hskip 6pt$}}

\def \C{\hbox{$C\hskip -5pt \vrule height 6pt depth 0pt \hskip 6pt$}}

\def\qed{\ \ \ifhmode\unskip\nobreak\fi\ifmmode\ifinner
         \else\hskip5pt\fi\fi
 \hbox{\hskip5pt\vrule width4pt height6pt depth1.5pt\hskip 1 pt}}
\def\a{\alpha}
\def\b{\beta}
\def\d{\delta}

\def\g{\gamma}
\def\G{\Gamma}
\def\l{\lambda}

\def\Vir{\mbox{Vir}}

\def\si{\sigma}
\def\sc{\scriptstyle}
\def\ssc{\scriptscriptstyle}
\def\F{\hbox{$I\hskip -4pt F$}}\def \Z{\hbox{$Z\hskip -5.2pt Z$}}
\def\dis{\displaystyle}
\def\cl{\centerline}
\def\DD{{\cal D}}

\def\nl{\newline}
\def\ol{\overline}

\def\wh{\widehat}
\def\rar{\rightarrow}

\def\bs{\backslash}

\def\rb{\raisebox}
\def\vs{\vspace*}
\def\ra{\rangle}
\def\la{\langle}
\def\ni{\noindent}
\def\hi{\hangindent}
\def\ha{\hangafter}
\def\WW{{\cal W}}
\def\AA{{\cal A}}

\begin{document}
\cl
{\large Classification of Quasifinite Modules over
        the Lie Algebras of Weyl Type}
\vs{5pt}\cl{(appeared in {\it Adv.~Math.} {\bf174} (2003),
57-68.)} \vs{10pt}\par \cl{{Yucai Su}\footnote{Supported by a NSF
grant 10171064 of China and a Fund from National Education
Ministry of China.\vs{4pt}\nl\hspace*{4.5ex}{\it Mathematics
Subject Classification (1991):} 17B10, 17B65, 17B66, 17B68}}
\cl{\small\it Department of Mathematics, Shanghai Jiaotong
University, Shanghai 200030, P.~R.~China} \vs{10pt}\par {\small
For a nondegenerate additive subgroup $\G$ of the $n$-dimensional
vector space $\F^n$ over an algebraically closed field $\F$ of
characteristic zero, there is an associative algebra and a Lie
algebra of Weyl type $\WW(\G,n)$ spanned by all differential
operators $uD_1^{m_1}\cdots D_n^{m_n}$ for $u\in\F[\G]$ (the group
algebra), and $m_1,...,m_n\ge0$, where $D_1, ...,D_n$ are degree
operators. In this paper, it is proved that an irreducible
quasifinite $\WW(\Z,1)$-module is either a highest or lowest
weight module or else a module of the intermediate series;
furthermore, a classification of uniformly bounded
$\WW(\Z,1)$-modules is completely given. It is also proved that an
irreducible quasifinite $\WW(\G,n)$-module is a module of the
intermediate series and a complete classification of quasifinite
$\WW(\G,n)$-modules is also given, if $\G$ is not isomorphic to
$\Z$. } \vs{10pt}\par \cl{1. \ INTRODUCTION}
\par
Let $n$ be a positive integer. A (classical)
{\it Weyl algebra of rank $n$} is the associative algebra
$A_n^+=\C[t_1,...,t_n,{\ptl\over\ptl{\ssc\,}t_1},...,
{\ptl\over\ptl{\ssc\,}t_n}]$ or
$A_n=\C[t_1^{\pm1},...,t_n^{\pm1},{\ptl\over\ptl{\ssc\,}t_1},...,
{\ptl\over\ptl{\ssc\,}t_n}]$ of differential operators over the complex
field $\C$.
The Lie algebra with $A_n$ as the underlined vector space and
the commutator as the Lie bracket is called a {\it Lie algebra of Weyl
type}, and denoted by $\WW(n)$.
The Lie algebra $\WW(n)$ is a
$\Z^n$-graded Lie
algebra $\WW(n)=\oplus_{\a\in\sZ^n}\WW(n)_\a$ with the grading space
$\WW(n)_\a$ spanned by
$$
t^\a D^\mu=t_1^{\a_1}\cdots t_n^{\a_n}D_1^{\mu_1}\cdots
D_n^{\mu_n}\mbox{ \ \ for \ \ }\mu=(\mu_1,...,\mu_n)\in\Z_+^n,
$$
where $\a=(\a_1,...,\a_n)\in\Z^n$ and $D_i=
t_i{\ptl\over\ptl{\ssc\,}t_i}$.
It is known [8,14,19] that $\WW(n)$ has a nontrivial universal
central extension if and only if $n=1$. The universal central extension
$\wh\WW(1)$ of $\WW(1)$ is the well-known Lie algebra $\WW_{1+\infty}$ of
the $\WW$-infinity algebras, which arise naturally in various physical theories
such as conformal field theory, the theory of the quantum Hall effect,
etc.~and which have received intensive studies in the literature
(see for example, [1-6,8-14,16,19-21]).
\par
Let $\F$ be an algebraically closed field of characteristic zero.
Consider the vector space $\F^n$. An element in $\F^n$ is denoted by
$\a=(\a_1,...,\a_n)$.
Let $\G$ be an additive subgroup of $\F^n$ such that $\G$ is
{\it nondegenerate}, i.e., it contains an $\F$-basis of $\F^n$.
Let $\F[\G]$ denote the group algebra of $\G$ spanned by $\{t^\a\,|\,
\a\in\G\}$ with the algebraic operation $t^\a\cdot t^\b=t^{\a+\b}$ for
$\a,\b\in\G$. We define the {\it degree operators} $D_i$
to be the derivations of $\F[\G]$ determined
by $D_i:t^\a\mapsto\a_it^\a$ for $\a\in\G$, where $i=1,...,n$.
The {\it Lie algebra of generalized Weyl type of rank $n$} (or simply
a {\it Lie algebra of Weyl type}, also called a {\it Lie algebra of
generalized differential operators}) is a
tensor product space $\WW(\G,n)=\F[\G,D_1,...,D_n]=
\F[\G]\otimes\F[D_1,...,D_n]$ of
the commutative associative algebra $\F[\G]$ with the polynomial algebra
$\F[D_1,...,D_n]$
of $D_1,...,D_n$, which is spanned by
$\{t^\a D^\mu\,|\,\a\in\G,\,\mu\in\Z_+^n\}$, where $D^\mu$ stands for
$\prod_{i=1}^nD_i^{\mu_i}$, with the Lie bracket:
$$
\begin{array}{c}
[t^\a D^\mu,t^\b D^\nu]=
(t^\a D^\mu)\cdot(t^\b D^\nu)
-(t^\b D^\nu)\cdot(t^\a D^\mu),\mbox{ \ \ and}
\vs{4pt}\\ \dis
(t^\a D^\mu)\cdot(t^\b D^\nu)
=\sum_{\l\in\sZ_+^n}\bigl(\!\begin{array}{c}\mu\\ \l\end{array}\!\bigr)
[\b]^\l t^{\a+\b}D^{\mu+\nu-\l},\mbox{ \ \ where}
\vs{4pt}\\ \dis
{[\b]}^\l=
\prod_{i=1}^n\b_i^{\l_i},\;\;\;\;
\bigl(\!\begin{array}{c}\mu\\ \l\end{array}\!\bigr)
=\prod_{i=1}^n
\bigl(\!\begin{array}{c}\mu_i\\ \l_i\end{array}\!\bigr),
\;\;\;\;
\bigl(\!\begin{array}{c}i\\ j\end{array}\!\bigr)
=\left\{\!
\begin{array}{ll}
{i(i-1)\cdots(i-j+1)\over j!}&\mbox{if \ }j\ge0,\vs{2pt}\\
0&\mbox{otherwise}.
\end{array}\right.
\end{array}
\eqno\begin{array}{r}(1.1)\!\!\vs{10pt}\\ (1.2)\!\!\vs{20pt}\\ \ \\ \
\end{array}\vs{-6pt}$$
The associative algebra with the underlined vector space $\WW(\G,n)$ and
the product (1.2) is called a
{\it generalized Weyl algebra of rank $n$}, denoted by $\AA(\G,n)$.
Then the classical Lie algebra $\WW(n)$ of Weyl type is simply
the Lie algebra $\WW(\Z^n,n)$ by our definition.
\par
Clearly $\WW(\G,n)$ is a $\G$-graded Lie algebra $\WW(\G,n)\!=\!
\oplus_{\a\in\G}\WW(\G,n)_\a$ with the grading space
$$
\WW(\G,n)_\a={\rm span}\{t^\a D^\mu\,|\,\mu\in\Z_+^n\}
\mbox{ \ \ for \ \ }\a\in\G.
\eqno(1.3)$$
It is proved in [19] that $\WW(\G,n)$ has a nontrivial
universal central extension if
and only if $n=1$. The universal central extension $\wh\WW(\G,1)$ of
$\WW(\G,1)$ is defined as follows: The Lie bracket (1.1) is replaced by
$$
[t^\a[D]_\mu,t^\b[D]_\nu]\!=\!
(t^\a[D]_\mu)\cdot(t^\b[D]_\nu)
\!-\!(t^\b[D]_\nu)\cdot(t^\a[D]_\mu)\!+\!\d_{\a,-\b}(-1)^\mu\mu!\nu!
\big(\!\begin{array}{c}\a\!+\!\mu\\ \mu\!+\!\nu\!+\!1\end{array}\!\bigr)C,
\eqno(1.4)$$
for $\a,\b\in\G\subset\F,\,\mu,\nu\in\Z_+,$
where $[D]_\mu{\sc\!}={\sc\!}D(D{\sc\!}-{\sc\!}1)
\cdots(D{\sc\!}-{\sc\!}\mu{\sc\!}+{\sc\!}1)$, and $C$ is
a central element
of $\wh\WW(\G,1)$ (the 2-cocycle of $\WW(\G,1)$
corresponding to (1.4) seems to appear first in [10]).
\par
A $\WW(\G,n)$-module (or an $\AA(\G,n)$-module) $V$ is called a
{\it quasifinite module}
[1,6,8,11,13,16] if
$V=\oplus_{\a\in\G}V_\a$ is $\G$-graded such that
$\WW(\G,n)_\a V_\b\subset V_{\a+\b}$ for $\a,\b\in\G$ and such that
each grading space $V_\a$ is finite-dimensional
(this is equivalent to saying that a quasifinite module is a module having
finite dimensional generalized weight spaces with respect to the commutative
subalgebra $\WW_0$).
A quasifinite module $V$ is called a {\it uniformly bounded module} if
the dimensions of grading spaces $V_\a$ are uniformly bounded, i.e.,
there
exists a positive integer $N$ such that ${\rm dim\ssc\,}V_\a\le N$ for all
$\a\in\G$. A quasifinite module $V$ is called a {\it trivial module}
if $\WW(\G,n)$ acts trivially on $V$.
\par
Clearly an $\AA(\G,n)$-module is a $\WW(\G,n)$-module, but not the converse.
Thus it suffices to consider $\WW(\G,n)$-modules.
The quasifinite highest weight modules over $\wh\WW(1)=\wh\WW(\Z,1)$
were intensively studied in [1,6,9-11,13,16,21]
(it is also worth mentioning that Block [3] studied arbitrary
irreducible modules over the classical Weyl algebra $A^+_1$ of rank 1).
\par
Let $p\ge1$.
Denote by $M_{p\times p}(\F)$ the set of $p\times p$ matrices with entries
in $\F$.
Denote
$$
M^n_{p\times p}(\F)=\{G=(G_1,...,G_n)\,|\,
G_i\in M_{p\times p}(\F),\,G_iG_j=G_jG_i\mbox{ \,for \,}i,j=1,...,n\},
$$
the set of $n$-tuples of commuting $p\times p$ matrices.
Denote by ${\bf1}_p$ the $p\times p$ identity matrix.
Denote $\one=({\bf1}_p,...,{\bf1}_p)\in M^n_{p\times p}(\F)$. Let
$G=(G_1,...,G_n)\in M^n_{p\times p}(\F)$. Then
one can define a quasifinite $\WW(\G,n)$-module $A_{p,G}$ as follows:
it has a basis
$\{y^{(i)}_\a\,|\,\a\in\G,\,i=1,...,p\}$ such that
$$
\begin{array}{c}
(t^\a D^\mu)Y_\b=Y_{\a+\b}[\b\cdot\one+G]^\mu
\mbox{ \ \ for \ \ }\a,\b\in\G,\,\mu\in\Z_+^n,\mbox{ \ \ where}
\vs{4pt}\\ \dis
Y_\b=(y^{(1)}_\b,...,y^{(p)}_\b),\;\;\;\;
[\b\cdot\one+G]^\mu=\prod_{i=1}^n(\b_i{\sc\!}\cdot{\sc\!}{\bf1}_p+G_i)^{\mu_i}.
\end{array}
\eqno\begin{array}{r}(1.5)\!\!\vs{16pt}\\ \ \end{array}\vs{-6pt}
$$
Here $\b_i{\sc\!}\cdot{\sc\!}{\bf1}_p$ denotes the scalar multiplication
of the identity matrix.
Clearly $A_{p,G}$ is also an $\AA(\G,n)$-module.\,By [20],\,there$\,
$exists\,a\,Lie algebra
isomorphism\,$\si{\sc\!}:\!\WW(\G,n)\!\cong\!\WW(\G,n)$\,such that
$$
\si(t^\a D^\mu)=(-1)^{|\mu|+1}([D^\mu,t^\a]+t^\a D^\mu)
\mbox{ \ \ for \ \ }\a\in\G,\,\mu\in\Z_+,
\eqno(1.6)$$
where $|\mu|=\sum_{i=1}^n\mu_p$ is the {\it level} of $\mu=(\mu_1,...,\mu_n)
\in\Z_+^n$. This isomorphism is uniquely
determined by $\si(t^\a)=-t^\a,\si(D^j_i)
=(-1)^{j+1}D_i$ for $\a\in\G,\,i=1,...,n,\,j\in\Z_+$.
Using this isomorphism, we have another $\WW(\G,n)$-module $\ol A_{p,G}$,
called the {\it twisted module of $A_{p,G}$}, for the
pair $(p,G)$, defined by
$$
(t^\a D^\mu)Y_\b=(-1)^{|\mu|+1}Y_{\a+\b}[(\a+\b)\cdot\one+G]^\mu
\mbox{ \ \ for \ \ }\a,\b\in\G,\,\mu\in\Z_+^n.
\eqno(1.7)$$
Obviously, $\ol A_{p,G}$ is not an $\AA(\G,n)$-module.
When $\G=\Z,n=p=1$, the above modules $A_{p,G},\ol A_{p,G}$ were
obtained in [21].
Clearly, $A_{p,G}$ or $\ol A_{p,G}$
is indecomposable if and only if at least one $G_i$ is
an indecomposable matrix
(here a $p\times p$ matrix $B$ is called
{\it indecomposable} if there does
not exist an invertible matrix $P$ such that $P^{-1}BP=
{\rm diag}(B_1,B_2)$ for some $p_i\times p_i$ matrices $B_i$
with $p_i<p,\,i=1,2$)
and $A_{p,G}$ or $\ol A_{p,G}$ is irreducible if and only if $p=1$.
When $p=1$, we refer $A_G=A_{1,G}$ or $\ol A_G=\ol A_{1,G}$ to as a
{\it module of the intermediate series}
(a notion borrowed from that of modules over the Virasoro algebra,
cf.~[15]).
\par
Since each grading space in (1.3) is infinite-dimensional, the
classification of quasifinite modules is thus a nontrivial problem, as
pointed in [11,13].
The aim of this paper is to prove the following theorem (the analogous
results for the affine Lie algebras and the Virasoro algebra were obtained
in [7,15]).
\par
{T$\sc\rm HEOREM$ 1.1.}\
{\it
(i) A uniformly bounded module over $\WW(1)=\WW(\Z,1)$ or over
$\wh\WW(1)=\wh\WW(\Z,1)$ is a direct sum of a trivial module, a
module $A_{p,G}$ and a module $\ol A_{p',G'}$ for some positive
integers $p,p'$ and some $G\in M^n_{p\times p}(\F),G'\in M_{p'\times p'}(\F)$
(in the central extension case,
the central element $C$ acts trivially on a uniformly bounded module);
a quasifinite irreducible module is
either a highest or lowest weight module or else
a module of the intermediate series.
\par
(ii) Suppose $\G$ is not isomorphic to $\Z$.
A quasifinite module
over $\WW(\G,n)$ or over $\wh\WW(\G,1)$ is
a direct sum of a trivial module and a
uniformly bounded module; a uniformly bounded module is
a direct
sum of a trivial module, a module $A_{p,G}$ and a module
$\ol A_{p',G'}$ for some positive integers $p,p'$ and some
$G\in M^n_{p\times p}(\F),G'\in M^n_{p'\times p'}(\F)$.
In particular, a nontrivial quasifinite irreducible module is a module of
the intermediate series.}
\par
Thus we in particular obtain that an indecomposable uniformly
bounded $\WW(\G,n)$-module is simply an $\AA(\G,n)$-module (if the
central element $t^0D^0$ acts by $1$) or its twist (if $t^0D^0$
acts by $-1$), and that there is an equivalence between the
category of uniformly bounded $\AA(\G,n)$-modules without the
trivial composition factor and the category of the finite
dimensional $\WW_0$-modules obtained by restriction to any nonzero
graded subspace.
\par
A {\it composition series of a module} $V$ is a finite or infinite
series of submodules $V=V^{(0)}\supset
V^{(1)}\supset V^{(2)}\supset...\supset \{0\}$ such that each $V^{(i)}/V^{(i+1)}$, called
a {\it composition factor}, is irreducible.
\par
{R$\sc\rm EMARK$ 1.2.}\
Note that the definition of quasifiniteness does not require that
$V$ is a {\it weight module} (i.e., the actions of $D_i,\,i=1,...,n$
on $V$ are diagonalizable). If we require
that $V$ is a weight module,
then each $G_i$ in (1.5) is diagonalizable, and thus all uniformly bounded
modules are completely reducible. Also note that in Theorem 1.1(ii),
if a module have infinite number of the trivial composition factor, then
it is not necessarily uniformly bounded since any $\G$-graded
vector space can be defined
as a trivial module.
\par
We shall prove Theorem 1.1(i) and (ii) in Sections 2 and 3 respectively.
\vs{10pt}\par
\cl{2. \ PROOF OF THEOREM 1.1(i)}
\par
For convenience, we shall only work on non-central extension case since
the proof of central extension case is exactly similar.
In this section we shall consider the Lie algebra
$\WW(\Z,1)={\rm span}\{t^iD^j\,|\,(i,j)\in\Z\times\Z_+\}$, which is
now denoted by $\WW$. Then $\WW$ is $\Z$-graded
$\WW=\oplus_{i\in\sZ}\WW_i$ with
$\WW_i={\rm span}\{t^iD^j\,|\,j\in\Z_+\}$, and
it has a triangular decomposition $\WW=\WW_+\oplus\WW_0\oplus\WW_-$
with $\WW_+=\oplus_{i>0}\WW_i,\WW_-=\oplus_{i<0}\WW_i$. Observe that
$\WW_+$ is generated by the adjoint action of $L_{1,0}$ on $\WW_0$ and
that $\ad_{L_{1,0}}$ is locally nilpotent on $\WW$,
where for convenience, we denote
$$
L_{i,j}=t^i[D]_j\mbox{ \ \ for \ \ }(i,j)\in\Z\times\Z_+.
$$
Note that in fact $L_{i,j}=t^{i+j}({d\over dt})^j$ (cf.~notation in (1.4)).
Suppose $V=\oplus_{i\in\sZ}V_i$ is a quasifinite module over $\WW$.
For any $a\in\F$, we denote
\vs{-4pt}$$
V(a)=\bigoplus_{i\in\sZ}V(a)_i,\mbox{ \ \ where \ \ }
V(a)_i=\{v\in V_i\,|\,L_{0,1}v=(a+i)v\}.
\vs{-6pt}$$
Since $[L_{0,1},L_{i,j}]=iL_{i,j}$, it is straightforward to verify that
$V(a)$ is
a submodule of $V$. Since $V_i$ is finite-dimensional, there exists at least
an eigenvector of $L_{0,1}$ in $V_i$, i.e., $V(a)\ne0$ for some $a\in\F$.
If $V$ is irreducible, then $V=V(a)$ for some $a\in\F$, i.e.,
$L_{0,1}$ is diagonalizable on $V$ and so $V$ is a weight module.
Since $L_{0,0}=1$ is a central element and $L_{0,0}|_{V_0}$ has at least an
eigenvector, we must have
$$
L_{0,0}|_V=c{\sc\!}\cdot{\sc\!}{\bf1}_V\mbox{ \ \ for some \ \ }c\in\F
\mbox{ \ \ if \ $V$ \ is indecomposable}.
\eqno(2.1)$$
\par
{P$\sc\rm ROPOSITION$ 2.1.}\
{\it
Suppose $V$ is a quasifinite irreducible $\WW$-module which
is neither a highest
nor a lowest weight module. Then
$L_{1,0}:V_i\rar V_{i+1}$ and $L_{-1,0}:V_i\rar V_{i-1}$ are injective and
thus bijective for all $i\in\Z$. In particular, $V$ is uniformly bounded.
}\par
{\it Proof.}\
Say $L_{1,0}v_0=0$ for some $0\ne v_0\in V_i$.
By shifting the grading index if necessary, we can suppose $i=0$.
Since
$L_{0,0}|_{V_0},L_{0,1}|_{V_0},...$ are linear transformations on
the finite-dimensional vector space $V_0$, there exists
$s\ge2$ such that for all $k\ge s$, $L_{0,k}|_{V_0}$ are linear
combinations of $L_{0,p}|_{V_0},\,0\le p<s$.
This implies that $\WW_0v_0=Sv_0$, where
$S={\rm span}\{L_{0,p}\,|\,0\le p<s\}.$
Recall that $\ad_{L_{1,0}}$ is locally nilpotent such that $\WW_k=
\ad_{L_{1,0}}^k(\WW_0)$ for $k>0$. Choose $m>0$ such that
$\ad_{L_{1,0}}^m(S)=0$, then for $k\ge m$, one has
$$
\WW_kv_0=(\ad_{L_{1,0}}^k(\WW_0))v_0=L_{1,0}^k\WW_0v_0=L_{1,0}^kSv_0
=(\ad_{L_{1,0}}^k(S))v_0=0.
$$
This means
$$
\WW_{[m,\infty)}v_0=0,
\eqno(2.2)$$
where in general,
for any $\Z$-graded space $N$, we use notations $N_+,N_-,N_0$ and
$N_{[p,q)}$
to denote the subspaces spanned by elements of degree $k$ with
$k>0,k<0,k=0$ and $p\le k< q$ respectively.
For any subspace $M$ of $\WW$, we use $U(M)$ to denote the subspace,
which is the span of standard monomials with respect to a basis of $M$,
of the universal enveloping algebra of $\WW$.
Since $\WW=\WW_{[1,m)}+\WW_0+\WW_-+\WW_{[m,\infty)}$,
using the PBW theorem and the irreducibility of $V$,
we have
$$
V=U(\WW)v_0=U(\WW_{[1,m)})U(\WW_0+\WW_-)U(\WW_{[m,\infty)})v_0
=U(\WW_{[1,m)})U(\WW_0+\WW_-)v_0.
\eqno(2.3)$$
\par
Note that $V_+$ is a $\WW_+$-module. Let $V'_+$ be the $\WW_+$-submodule of
$V_+$ generated by $V_{[0,m)}$.
We want to
prove
$$
V_+=V'_+.
\eqno(2.4)$$
Let $x\in V_+$ have degree $k$.
If $0\le k<m$, then by definition $x\in V'_+$. Suppose $k\ge m$.
Using (2.3), $x$ is a linear combination of the form $u_1x_1$ with
$u_1\in\WW_{[1,m)},x_1\in V$.
Thus the degree ${\rm deg\ssc\,}u_1$ of $u_1$ satisfies
$1\le{\rm deg\ssc\,}u_1<m$,
so $0<{\rm deg\ssc\,}x_1=k-{\rm deg\ssc\,}u_1<k$.
By inductive hypothesis, $x_1\in V'_+$, and thus $x\in V'_+$. This proves
(2.4).
\par
(2.4) means that $V_+$ is finitely generated as a $\WW_+$-module.
Choose a basis $B$ of $V_{[0,m)}$, then for any $x\in B$, we have
$x=u_xv_0$ for some $u_x\in U(\WW)$.
Regarding $u_x$ as a polynomial with respect to a basis of $\WW$, by
induction on the polynomial degree and using the formula $[u,w_1w_2]
=[u,w_1]w_2+w_1[u,w_2]$ for $u\in\WW,\,w_1,w_2\in U(\WW)$,
we see that there exists
a positive integer $k_x$ large enough such that $k_x>m$ and
$[\WW_{[k_x,\infty)},u_x]\subset U(\WW)\WW_{[m,\infty)}$. Then
by (2.2), $\WW_{[k_x,\infty)}x
=[\WW_{[k_x,\infty)},u_x]v_0+u_x\WW_{[k_x,\infty)}v_0
=0$. Take $k={\rm max}\{k_x\,|\,x\in B\}$, then
$$
\WW_{[k,\infty)}V_+=\WW_{[k,\infty)}U(\WW_+)V_{[0,m)}=
U(\WW_+)\WW_{[k,\infty)}V_{[0,m)}=0.
$$
Since we have $\WW_+\subset\WW_{[k,\infty)}+[\WW_{[-k',0)},\WW_{[k,\infty)}]$
for some $k'>k$, we get $\WW_+V_{[k',\infty)}=0$. Now if $x\in
V_{[k'+m,\infty)}$, by (2.3), it is a sum of elements of the form $u_1x_1$
such that $u_1\in\WW_{[1,m)}$ and so $x_1\in V_{[k',\infty)}$, and
thus $u_1x_1=0$. This proves that $V$ has no degree $\ge k'+m$.
\par
Now let $p$ be maximal integer such that $V_p\ne0$.
Since $\WW_0$ is commutative, there exists a common eigenvector $v'_0$
of $\WW_0$ in $V_p$.
Then $v'_0$ is a highest weight vector of $\WW$, this contradicts
the assumption of the proposition.
\qed\par
{P$\sc\rm ROPOSITION$ 2.2.}\
{\it Suppose $V$ is an indecomposable
uniformly bounded $\WW$-module without the trivial
composition factor. Then $V$ is a module of the form
$A_{p,G}$ or $\ol A_{p,G}$.}
\par
{\it Proof.}\
First we claim that $L_{i,0}$ acts nondegenerately on $V$ for all $i\ne0$.
Say $L_{i_0,0}v_0=0$ for some $i_0>0$ and some $0\ne v_0\in V_0$. Denote
$\WW'=\WW(\Z i_0,1)={\rm span}\{L_{i_0j,k}\,|\,j,k\in\Z,k\ge0\}$, a
subalgebra
of $\WW$, which is clearly isomorphic to $\WW$. Let $V'$ be the
$\WW'$-submodule of $V$ generated by $v_0$.
Then by replacing $\WW$ and $V$ by $\WW'$ and $V'$ respectively in
the proof of Proposition 2.1, we see that $V'$ has a highest weight
vector $v'$ with respect to $\WW'$
(cf.~(2.3), the left-hand side of (2.3) is now $V'$).
Then $v'$ generates a highest weight $\WW'$-submodule $V''$.
Since a nontrivial highest weight $\WW'$-module is not uniformly bounded
(see for example [11]), $V''$ must be a trivial $\WW'$-module $V''=\F v'$.
Denote $\Vir={\rm span}\{L_{j,1}\,|\,j\in\Z\}$, a subalgebra of $\WW$,
which is isomorphic to the centerless Virasoro algebra, then $v'$ generates
a uniformly bounded weight $\Vir$-submodule $U$ of $V$
(note that in general we do not assume $V$ is a weight module over $\WW$).
Since $L_{i_0j,1}v'=0$ for
all $j\in\Z$, by a well-known result of a uniformly bounded weight
$\Vir$-module
(see for example, [15,18]), we have $L_{i,1}v'=0$ for all $i\in\Z$. But
$\WW$ is generated by $\WW'$ and $\Vir$, we obtain that $\F v'$ is a trivial
$\WW$-module. This is a contradiction with the assumption of the proposition.
This proves the claim.
\par
So, there exists $p\ge1$ such
that ${\rm dim\ssc\,}V_k=p$ for all $k\in\Z$, and one can choose a basis
$Y_0=(y^{(1)}_0,...,y^{(p)}_0)$ of $V_0$ and define a basis
$Y_k=(y^{(1)}_k,...,y^{(p)}_k)$ of $V_k$ by induction on $|k|$
such that
$$
L_{1,0}Y_k=Y_{k+1}\mbox{ \ \ for \ \ }k\in\Z.
\eqno(2.5)$$
Furthermore, by induction on $p$, we see that $V$ has a finite number of
composition factors.
\par
First note that
\vs{-4pt}$$
[t^i(\mbox{$d\over dt$})^j,t^k(\mbox{$d\over dt$})^\ell]=
\sum_{s}
\rb{-2pt}{\mbox{\Large(}}
\rb{-1pt}{\mbox{\large(}}\!\begin{array}{c}j\\ s
\end{array}\!\rb{-1pt}{\mbox{\large)}}[k]_s
-\rb{-1pt}{\mbox{\large(}}\!\begin{array}{c}\ell\\ s
\end{array}\!\rb{-1pt}{\mbox{\large)}}[i]_s
\rb{-2pt}{\mbox{\Large)}}
t^{i+k-s}(\mbox{$d\over dt$})^{j+\ell-s},
\vs{-4pt}\eqno(2.6)$$
where $[k]_j=k(k-1)\cdots(k-j+1)$ is a similar notation to $[D]_j$.
Now assume that $L_{i-j,j}Y_n=(t^i({d\over dt})^j)Y_n=Y_{n+i-j}P_{i,j,n}$
for some $P_{i,j,n}\in M_{p\times p}$.
In the following discussion, we remind the reader that
$t^i({d\over dt})^j$ is in the grading space $\WW_{i-j}$, not in
$\WW_i$.
Applying $[t^i({d\over dt})^j,t]
=jt^i({d\over dt})^{j-1}$ to
$Y_n$, we obtain $P_{i,j,n+1}-P_{i,j,n}=jP_{i,j-1,n}$. Thus induction on
$j$ gives
$$
\begin{array}{c}
P_{i,0,n}=P_i,\;\;\;P_{i,1,n}=\bar nP_i+Q_i,
\vs{4pt}\\
P_{i,2,n}=[\bar n]_2P_i+2\bar nQ_i+R_i,\;\;\;
P_{i,3,n}=[\bar n]_3P_i+3[\bar n]_2Q_i+3\bar nR_i+S_i,
\end{array}
\eqno(2.7)$$
for some $P_i,Q_i,R_i,S_i\in M_{p\times p}$,
where $\bar n=n+G$ for some fixed $G\in M_{p\times p}$
(here and below,
when the context is clear, we identify a scalar with the corresponding
$p\times p$ scalar matrix),
$[\bar n]_j$ is again a similar notation to $[D]_j$, and
$Q_1=0$ (we use notation $\bar n=n+G$ in order to have $Q_1=0$;
note from $[L_{0,1},L_{i,j}]=iL_{i,j}$ that $G$ commutes with all other
matrices involved in the following discussion).
Then by (2.1) and (2.5), $P_0=c,\,P_1=1$
and $P_i$ is invertible for $i\ne0$ by the
proof above. Applying $[{d\over dt},t^2]=2t$ to $Y_n$ and comparing
the coefficients of $\bar n^0$, we obtain
$2cP_2+[Q_0,P_2]=2$ (where $[P_2,Q_0]=P_2Q_0-Q_0P_2$ denotes the usual
commutator of matrices). Comparing the traces of matrices
in this equation shows $c\ne0$.
Thus all $P_i$ are invertible.
\par
Now we encounter the difficulty that though
$P_i,Q_i,R_i,S_i$ satisfy lots of relations, most nontrivial
relations are too complicated to be used; since the products of matrices
are not commutative nor the cancellation law holds in general,
there is still a problem in finding the solutions for $P_i,Q_i,R_i,S_i$.
Our strategy is first to find some relatively simple nontrivial
relations among $P_i$ (cf.~(2.9) below).
\par
First from $[t^i,t^j]=0$, we obtain that $[P_i,P_j]=0$.
Applying $[{d\over dt},t^i]=it^{i-1}$ to $Y_n$ and comparing
the coefficients of $\bar n^0$, we obtain that
$[Q_0,P_i]=-ic{\ssc\,}P_i+iP_{i-1}$.
By induction on $i\ge0$, we obtain an important
fact that
$P'_i=\sum_{s=0}^i(^i_s)(-c)^sP_s$, if not zero,
is an eigenvector for ${\rm ad}_{Q_0}$
with eigenvalue $-ic$. Since the operator ${\rm ad}_{Q_0}$ acting
on the finite
dimensional vector space $M_{p\times p}$ has only a finite number of
eigenvalues, we have $P'_i=0$ for $i>>0$. Let $q\ge0$
be the least number such
that
\vs{-4pt}$$
P'_i=\sum_{s=0}^i
\mbox{\large(}\!
\begin{array}{c}
i\\ s\end{array}\!\mbox{\large)}
(-c)^sP_s=0,
\vs{-4pt}\eqno(2.8)$$
for $i>q$. Note that for any $j\ge0$, we have
\vs{-4pt}$$
\begin{array}{ll}\dis
P''_{i,j}
\!\!\!\!&\dis
:=
\sum_{k=0}^j
(-1)^k
\mbox{\large(}\!
\begin{array}{c}
i\\ k\end{array}\!\mbox{\large)}
P'_k
=\sum_{k=0}^j
\sum_{s=k}^s
(-1)^k
\mbox{\large(}\!
\begin{array}{c}
i\\ k\end{array}\!\mbox{\large)}
\mbox{\large(}\!
\begin{array}{c}
k\\ s\end{array}\!\mbox{\large)}
(-c)^sP_s
\vs{4pt}\\ &\dis
=\sum_{s=0}^j
c^s
\mbox{\large(}\!
\begin{array}{c}
i\\ s\end{array}\!\mbox{\large)}
\mbox{\Large(}
\sum_{k=s}^j
\mbox{\large(}\!
\begin{array}{c}
i\!-\!s\\ k\!-\!s\end{array}\!\mbox{\large)}
(-1)^{k-s}
\mbox{\Large)}P_s
=\sum_{s=0}^j c^s
\mbox{\large(}\!
\begin{array}{c}
i\\ s\end{array}\!\mbox{\large)}
\mbox{\large(}\!
\begin{array}{c}
i\!-\!1\!-\!s\\ j\!-\!s\end{array}\!\mbox{\large)}
(-1)^{j-s}P_s.
\end{array}
\vs{-4pt}$$
In particular letting $j=i$, we obtain that $P''_{i,i}=c^iP_i$ and
letting $j=q$, by (2.8), we obtain
\vs{-5pt}$$
P_i=c^{-i}P''_{i,i}=c^{-i}P''_{i,q}=
c^{-i}\sum_{s=0}^q
\mbox{\large(}\!
\begin{array}{c}
i\\ s\end{array}\!\mbox{\large)}
\mbox{\large(}\!
\begin{array}{c}
i\!-\!1\!-\!s\\ q-s\end{array}\!\mbox{\large)}
(-1)^{q-s}c^sP_s,
\vs{-4pt}\eqno(2.9)$$
for $i>q$. If $0\le i\le q$, (2.9) holds trivially.
Observe that $(^i_s)(^{i-1-s}_{\ q-s})$ is a polynomial on $i$ of degree $q$
with the coefficient of $i^q$ being ${1\over s!(q-s)!}={1\over q!}(^q_s)$.
Thus $(-1)^qc^iP_i$ is a polynomial on $i$
(with coefficients in $M_{p\times p}$) of degree $q$ with the
coefficient of $i^q$ being
${1\over q!}\sum_{s=0}^q(^q_s)(-c)^sP_s$.
\par
By (2.6), we have
$$
3[t(\mbox{$d\over dt$})^2,[t(\mbox{$d\over dt$})^2,t^i]]
-2(2i\!-\!1)[t(\mbox{$d\over dt$})^3,t^i]+[i\!+\!1]_4t^{i-2}=0.
$$
Applying
this to $Y_n$ and comparing the coefficients of $\bar n^0$, we obtain
$$
\begin{array}{ll}
f(i):=\!\!\!\!&3(
\mbox{\small(}[i\!-\!1]_2\!+\!R_1\mbox{\small)}
\mbox{\small(}[i]_2\!+\!R_1\mbox{\small)}P_i
-2\mbox{\small(}[i\!-\!1]_2\!+\!R_1\mbox{\small)}P_iR_1
\!+\!P_i\mbox{\small(}2\!+\!R_1\mbox{\small)}R_1)
\vs{5pt}\\ &
-2(2i\!-\!1)(
\mbox{\small(}[i]_3\!+\!3iR_1+S_1
\mbox{\small)}P_i
\!-\!P_iS_1)
\!+\![i\!+\!1]_4P_{i-2}\!=\!0.
\end{array}
\eqno(2.10)$$
Using (2.9) in (2.10), we obtain that
$(-1)^qc^if(i)$ is a polynomial on $i$ of degree at most $q+4$.
By comparing the coefficients of $i^{q+4}$, we obtain
\vs{-8pt}$$
{1\over q!}(c^2-1)\sum_{s=0}^q
\mbox{\large(}\!\begin{array}{c}q\\ s\end{array}\!\mbox{\large)}
(-c)^sP_s=0.
\vs{-8pt}\eqno(2.11)$$
If necessary, by using the isomorphism in (1.6)
(which interchanges $A_{p,G}$ with $\ol A_{p,G}$), we can
always suppose $c\ne-1$.
\par
Assume that $c\ne1$. Then (2.11) shows
that (2.8) also holds for $i=q$. Thus the minimality of $q$
implies that $q=0$ and then (2.9) implies
\vs{-5pt}$$
P_i=c^{1-i},
\vs{-5pt}\eqno(2.12)$$
for $i\ge0$.
Assume that $c=1$. Using (2.9) in (2.10) and comparing the coefficients
of $i^{q+3}$, we again obtain that (2.8) holds for $i=q$. Thus in any case,
$q=0$ and we have (2.12).
\par
For any $i\in\Z$, choose $j>0$ such that $i+j-1>0$.
Applying $[t^i{d\over dt},t^j]=jt^{i+j-1}$ to $Y_n$, we obtain
$-j(P_jP_i-P_{i+j-1})=[Q_i,P_j]=0$. This implies that (2.12)
holds for all $i\in\Z$.
\par
Now using (2.12) and
applying
$[({d\over dt})^2,t^{i+1}]=2(i\!+\!1)t^i{d\over dt}
+[i\!+\!1]_2t^{i-1}$ to $Y_n$, by comparing
the coefficients of $\bar n^0$, we obtain
$$
(i+1)(\mbox{\small(}iP_0+2Q_0\mbox{\small)}P_{i+1}
-\mbox{\small(}2Q_i+iP_{i-1}\mbox{\small)})
=[P_{i+1},R_0]=0.
$$
Thus
$Q_i=\mbox{$1\over2$}c^{1-i}(1-c)i+c^{-i}Q_0$ for $i\ne-1.$
Letting $i=1$, since $Q_1=0$, we obtain that $Q_0=-{1\over2}c(1-c)$.
Hence
\vs{-5pt}$$
Q_i=\mbox{$1\over2$}c^{1-i}(1-c)(i-1)\mbox{ \ for \ }i\ne-1.
\vs{-5pt}\eqno(2.13)$$
Applying
$[({d\over dt})^3,t^{i+1}]=3(i\!+\!1)t^i({d\over dt})^2
+3[i\!+\!1]_2t^{i-1}{d\over dt}+[i\!+\!1]_3t^{i-2}$ to $Y_n$, we obtain
$$
(i+1)(\mbox{\small(}
[i]_2P_0+3iQ_0+3R_0
\mbox{\small)}P_{i+1}
-\mbox{\small(}
3R_i+3iQ_{i-1}+[i]_2P_{i-2}\mbox{\small)})
=[P_{i+1},S_0]=0.
$$
Thus
$$
R_i=\mbox{$1\over6$}c^{1-i}(c-1)(
5-4c+\mbox{\small(}c-2\mbox{\small)}i)i+c^{-i}R_0\mbox{ \,for \,}i\ne-1.
\eqno(2.14)$$
Using
$[({d\over dt})^2,t^{i+1}{d\over dt}]=2(i\!+\!1)t^i({d\over dt})^2
+[i\!+\!1]_2t^{i-1}{d\over dt},$
we obtain
$$
([i]_2P_0+2iQ_0)Q_i+2P_{i+1}R_0-(2\mbox{\small(}i\!+\!1\mbox{\small)}R_i
+[i\!+\!1]_2Q_{i-1})=[Q_{i+1},R_0]=0.
\eqno(2.15)$$
Using (2.12)-(2.14) in (2.15), by comparing the coefficients of $i^3,i^2,i$,
we obtain $c=1,R_0=0$. Then (2.12)-(2.15) show that $P_i=1,Q_i=R_i=0$.
\par
Thus (2.7) shows
that $(t^i({d\over dt})^j)Y_n=Y_{n+i-j}[n+G]_j$
for $j\le 2$. Equivalently,
$(t^iD^j)Y_n=Y_{n+i}(n+G)^j$ for $j\le2$.
Since $\WW$ is generated by $\{t^iD^j\,|\,i\in\Z,0\le j\le2\}$, we obtain
that $V$ is the module $A_{p,G}$ defined in (1.5) (if we have used the
isomorphism in (1.6) in the above proof,
then $V$ is the module $\ol A_{p,G}$).
\qed\par
{\it Proof of Theorem 1.1(i).} \,
Assume that $V$ is a nontrivial indecomposable uniformly bounded module
over $\WW$. Then $V$ has at least a nontrivial composition factor, so
Proposition 2.2 means that $L_{0,0}$ must acts as a nonzero scalar on $V$.
Thus $V$ can not contain a trivial composition factor. Thus Theorem 1.1(i)
follows from Propositions 2.1 and 2.2.
\qed\vs{10pt}\par
\cl{3.\ PROOF OF THEOREM 1.2(ii)}
\par
In this section, we set $\WW=\WW(\G,n)$, where $\G$ is not isomorphic to
$\Z$.
Denote ${\sl Witt}={\rm span}\{t^\a D_i\,|\,\a\in\G,i=1,...,n\}$.
Then ${\sl Witt}$ is a Witt algebra
of rank $n$. In particular, if $n=1$, ${\sl Witt}$ is a (generalized
centerless) Virasoro algebra.
\par
Denote $\DD={\rm span}\{D_i\,|\,i=1,...,n\}$.
We can define an inner product on $\G\times\DD$ by
\vs{-5pt}$$
\la\a,d\ra=\sum_{i=1}^n\a_id_i\mbox{ \ \ for \ \ }
\a=(\a_1,...,\a_n)\in\G,\,d=\sum_{i=1}^nd_iD_i\in\DD.
\vs{-5pt}\eqno(3.1)$$
Then $\la\cdot,\cdot\ra$ is {\it nondegenerate} in the sense that
if $\la\a,\DD\ra=0$ for some $\a\in\G$
then $\a=0$ and if $\la\G,d\ra=0$ for some $d\in\DD$ then $d=0$.
\par
{P$\sc\rm ROPOSITION$ 3.1.}\
{\it
Any quasifinite $\WW$-module $V'$ with a finite number of the trivial
composition factor is a uniformly bounded module.}
\par
{\it Proof.}\
Let $V$ be a nontrivial composition factor of $V'$.
Then $V$ is a weight module over $\WW$ (i.e., $D_i$ act
diagonalizable on $V$ for $i=1,...,n$).
Regarding $V$ as a module over ${\sl Witt}$, then $V$ is a quasifinite weight
module over ${\sl Witt}$. If there exists some group embedding
$\Z\times\Z\rar\G$,
then by [18], $V$ is uniformly bounded.
If there does not exist a
group embedding $\Z\times\Z\rar\G$, then $n=1$ and $\G$ is a rank one
group with infinite generators (since $\G\not\cong\Z$).
By choosing a total ordering on $\G$ compatible with its group structure
one can prove (as in the proof of Proposition 2.1 or using similar
arguments as in [17] since in this case the group $\G$ behaves just
like the additive group $\Q$ of the rational numbers) that $V$ is
uniformly bounded. Hence in any case, $V$ is uniformly bounded.
Thus $V'$ has only a finite number of
nontrivial composition factors and so it is uniformly bounded.
\qed\par
For any pair $(\a,d)$ with $\a\in\G,\,d\in\DD$ such that $\la\a,d\ra\ne0$,
we have a Lie algebra of Weyl type, denoted by $\WW(\a,d)$, spanned by
$\{t^{i\a}d^j\,|\,i\in\Z,j\in\Z_+\}$,
which is isomorphic to $\WW(\Z,1)$.
\par
{P$\sc\rm ROPOSITION$ 3.2.}\
{\it
Let $V$ be a uniformly bounded $\WW$-module without the trivial composition
factor. Then $t^\a\cdot v\ne0$ for all $\a\in\G\bs\{0\},\,v\in V\bs\{0\}$.
}
\par
{\it Proof.}\
Suppose $t^\b\cdot v_0=0$ for some $\b\ne0,v_0\ne0$. By shifting the grading
index if necessary, we can suppose $v_0$ has degree $0$.
Let $V'_0=\{v\in V_0\,|\,t^\b\cdot v=0\}$. Then $V'_0$ is invariant under the
action of $\DD$ since $[\DD,t^\b]\subset\F t^\b$. Thus we can find a common
eigenvector (i.e., a weight vector), denoted again by $v_0$, of $\DD$ in
$V'_0$.
Thus $t^\b\cdot v_0=0$ and $\DD\cdot v_0\subset\F v_0$.
For any $d'\in\DD$ with $\la\b,d'\ra\ne0$, considering the
$\WW(\b,d')$-submodule $V''$ of $V$ generated by $v_0$, by Theorem 1.1(i),
and by (1.5) and (1.7), $V''$ must be the trivial submodule $\F v_0$.
In particular, $d'\cdot v_0=0$ for all
$d'$ with $\la\b,d'\ra\ne0$. Since such $d'$ span $\DD$, we obtain
$\DD\cdot v_0=0$.
If $t^\a\cdot v_0=0$ for all
$\a\in\G\bs\{0\}$, then by Theorem 1.1(i), $\F v_0$ is a trivial
submodule over $\WW(\a,d)$ for all pairs $(\a,d)$
with $\la\a,d\ra\ne0$, and so
it must be a trivial submodule over $\WW$. Thus $t^\g\cdot v_0\ne0$ for some
$\g\in\G\bs\{0\}$. Choose $d\in\DD$ such that
$\la\b,d\ra\ne0\ne\la\g,d\ra$.
Then by Theorem 1.1(i) and (1.5), $v_0$ generates a
submodule of the intermediate series over $\WW(\g,d)$ such that
$(t^{i\g} d^j)\cdot v_0=0$ for $j\ne0$ and all $i\in\Z$
(note that the coefficient of the right-hand side of (1.5) does not
depend on $\a$).
So $t^{i\b+j\g}\cdot v_0=
\la i\b,d\ra^{-1}[t^{j\g} d,t^{i\b}]\cdot v_0=0$ if $i\ne0$. Similarly,
$(t^{i\b+j\g}d)\cdot v_0=0$ if $i\ne0$.
Choose some $i,j$ with $i\ne0$ such that $a=\la i\b+j\g,d\ra\ne0$,
then $t^\g\cdot v_0=a^{-1}[t^{-i\b+(1-j)\g}d,t^{i\b+j\g}]\cdot v_0=0$,
a contradiction.
\qed\par
{\it Proof of Theorem 1.1(ii).}\
Suppose $V$ is an indecomposable
uniformly bounded module without the trivial composition
factor. Choose a basis $Y_0=(y_0^{(1)},...,y_0^{(p)})$ of $V_0$ and
by Proposition 3.2, we can define
basis $Y_\a=t^a\cdot Y_0$ of $V_\a$ for all $\a\in\G\bs\{0\}$.
For any $\b\in\G\bs\{0\}$ with $\la\b,D_i\ra\ne0$ for $i=1,...,n$,
set $V[\b]=\oplus_{k\in\sZ}V_{k\b}$, then
$V[\b]$ is a submodule of $V$ over $\WW(\b,D_i)$ of type
$A_{p,G_i}$ for some $G_i\in M_{p\times p}$, $i=1,...,n$.
Since $[\DD,\DD]=0$, we have $G_iG_j=G_jG_i$ for $i,j=1,...,n$.
Since $\{\b\in\G\,|\,\la\b,D_i\ra\ne0,\,i=1,...,n\}$
generates $\G$, it is straightforward to see that $V$ is a $\WW$-module
of type $A_{p,G}$ or $\ol A_{p,G}$
with $G=(G_1,...,G_n)\in M^n_{p\times p}(\F)$.
If $V$ contains a trivial composition factor, then $V$ must be a trivial
module as in the proof of Theorem 1.1(i).
Thus we obtain Theorem 1.1(ii).
\qed\vs{10pt}\par
\cl{ACKNOWLEDGEMENT}
\par
The author would like to thank Professor Kaiming Zhao for discussions.
\vs{10pt}
\par
\cl{REFERENCES}
\lineskip=5.3pt
\par\ni\hi3.5ex\ha1
[1] H.~Awata, M.~Fukuma, Y.~Matsuo and S.~Odake, Character and
    determinant formulae of quasifinite representations of the
    $\WW_{1+\infty}$ algebra,
    {\it Comm.~Math.~Phys.} {\bf172} (1995), 377-400.
\parskip=0.08truein
\par\ni\hi3.5ex\ha1
[2] S.~Bloch, Zeta values and differential operators on the circle,
    {\it J.~Alg.} {\bf182} (1996), 476-500.
\par\ni\hi3.5ex\ha1
[3] R.~E.~Block, The irreducible representations of the Lie algebra
    $sl(2)$ and of the Weyl algebra, {\it Adv.~Math.} {\bf39} (1981), 69-110.
\par\ni\hi3.5ex\ha1
[4] R.~Blumenhagen, W.~Eholzer, A.~Honecker, K.~Hornfeck and
    R.~H\"ubel, Unifying $\WW$-algebras,
    {\it Phys.~Lett.~B.} {\bf332} (1994), 51-60.
\par\ni\hi3.5ex\ha1
[5] P.~Bouwknegt and K.~Schoutens, $\WW$-symmetry in conformal
    field theory, {\it Phys.~Rep.} {\bf223} (1993), 183-276.
\par\ni\hi3.5ex\ha1
[6] C.~Boyallian, V.~Kac, J.~Liberati and C.~Yan,
    Quasifinite highest weight modules of the Lie algebra of matrix
    differential operators on the circle,
    {\it J.~Math.~Phys.} {\bf39} (1998), 2910-2928.
\par\ni\hi3.5ex\ha1
[7] V.~Chari, Integrable representations of affine Lie algebras,
    {\it Invent.~Math.} {\bf85} (1986), 317-335.
\par\ni\hi3.5ex\ha1
[8] B.~L.~Feigin, The Lie algebra $gl(\l)$ and the cohomology of
    the Lie algebra of differential operators, {\it Uspechi Math.~Nauk}
    {\bf35} (1988), 157-158.
\par\ni\hi3.5ex\ha1
[9] E.~Frenkel, V.~Kac, R.~Radul and W.~Wang,
    $\WW_{1+\infty}$ and $\WW(gl_N)$ with central charge $N$,
    {\it Comm.~Math.~Phys.} {\bf170} (1995), 337-357.
\par\ni\hi4ex\ha1
[10] V.~G.~Kac and D.~H.~Peterson, Spin and wedge representations of infinite
    dimensional Lie algebras and groups, {\it Proc.~Nat.~Acad.~Sci.~U.~S.~A.}
    {\bf 78} (1981), 3308-3312.
\par\ni\hi4ex\ha1
[11] V.~G.~Kac and A.~Radul, Quasi-finite highest weight modules over the
    Lie algebra of differential operators on the circle,
    {\it Comm.~Math.~Phys.} {\bf157} (1993), 429-457.
\par\ni\hi4ex\ha1
[12] V.~G.~Kac and A.~Radul, Representation theory of the vertex algebra
     $\WW_{1+\infty}$, {\it Trans. Groups} {\bf1} (1996), 41-70.
\par\ni\hi4ex\ha1
[13] V.~G.~Kac, W.~Wang and C.~H.~Yan, Quasifinite  representations of classical
     Lie subalgebras of $\WW_{1+\infty}$, {\it Adv.~Math.}
     {\bf139} (1998), 46-140.
\par\ni\hi4ex\ha1
[14] W.~Li, 2-Cocycles on the algebra of differential operators,
     {\it J.~Alg.} {\bf122} (1989), 64-80.
\par\ni\hi4ex\ha1
[15] O.~Mathieu, Classification of Harish-Chandra modules over the
     Virasoro Lie algebra, {\it Invent.~Math.} {\bf 107} (1992), 225-234.
\par\ni\hi4ex\ha1
[16] Y.~Matsuo, Free fields and quasi-finite representations of
     $\WW_{1+\infty}$, {\it Phys.~Lett.~B.} {\bf326} (1994), 95-100.
\par\ni\hi4ex\ha1
[17] V.~Mazorchuk, Classification of simple Harish-Chandra modules over
     $\Q$-Virasoro algebra, {\it Math.~Nachr.} {\bf209} (2000), 171-177.
\par\ni\hi4ex\ha1
[18] Y.~Su, Classification of Harish-Chandra modules over the
     higher rank Virasoro algebras, to appear.
\par\ni\hi4ex\ha1
[19] Y.~Su, 2-Cocycles on the Lie algebras of generalized differential
     operators, {\it Comm.~Alg.}, {\bf30} (2002), 763-782.
\par\ni\hi4ex\ha1
[20] Y.~Su and K.~Zhao, Isomorphism classes and automorphism groups of
     algebras of Weyl type, {\it Science in China A}, {\bf45}
     (2002), 953-963.
\par\ni\hi4ex\ha1
[21] K.~Zhao, The classification of a kind of irreducible Harish-Chandra
     modules over algebras of differential operators,
    {\it (Chinese) Acta Math.~Sinica} {\bf37} (1994), 332-337.
\end{document}